\documentclass[journal]{IEEEtran}

\usepackage{cite}

\usepackage{graphicx}
\usepackage{epstopdf}
\usepackage{pgf}
\usepackage{tikz}
\usepackage[utf8]{inputenc}
\usepackage{pgfplots}
\pgfplotsset{compat=1.14}
\usetikzlibrary{calc}
\usetikzlibrary{shapes.misc}
\usetikzlibrary{shapes.arrows}
\usetikzlibrary{shapes.symbols}
\usetikzlibrary{positioning}
\usetikzlibrary{arrows.meta}
\usetikzlibrary{shapes.geometric, arrows}

\usetikzlibrary{fit}
\graphicspath{{Images/}}
\DeclareGraphicsExtensions{.pdf,.jpeg,.png,.eps}
\usepackage[caption=false,font=footnotesize]{subfig}

\usepackage{amsmath}
\usepackage{amssymb}
\DeclareMathOperator{\E}{E}

\usepackage[short]{optidef}
\usepackage{bm}

\usepackage{siunitx}
\usepackage{eurosym}
\DeclareSIUnit{\EUR}{\text{\euro}}

\usepackage{algorithmicx}
\usepackage[noend]{algpseudocode}
\usepackage{algorithm}
\algrenewcommand\algorithmicindent{1.0em}%
\makeatletter
\newcommand{\algmargin}{\the\ALG@thistlm}
\makeatother

\algnewcommand\algorithmicparameter{\textbf{Parameters:}}
\algnewcommand\Parameter{\item[\algorithmicparameter]}
\algnewcommand{\parState}[1]{\State%
	\parbox[t]{\dimexpr\linewidth-\algmargin}{\strut #1\strut}}

\algnewcommand{\Initialize}[1]{%
	\State \textbf{Initialize:}
	\Statex \hspace*{\algorithmicindent}\parbox[t]{.8\linewidth}{\raggedright #1}
}

\usepackage{array}
\newcolumntype{L}[1]{>{\raggedright\let\newline\\\arraybackslash\hspace{0pt}}m{#1}}
\newcolumntype{C}[1]{>{\centering\let\newline\\\arraybackslash\hspace{0pt}}m{#1}}
\newcolumntype{R}[1]{>{\raggedleft\let\newline\\\arraybackslash\hspace{0pt}}m{#1}}

\usepackage{booktabs}
\usepackage{url}
\usepackage[hidelinks]{hyperref}
\expandafter\def\expandafter\UrlBreaks\expandafter{\UrlBreaks
	\do\a\do\b\do\c\do\d\do\e\do\f\do\g\do\h\do\i\do\j%
	\do\k\do\l\do\m\do\n\do\o\do\p\do\q\do\r\do\s\do\t%
	\do\u\do\v\do\w\do\x\do\y\do\z\do\A\do\B\do\C\do\D%
	\do\E\do\F\do\G\do\H\do\I\do\J\do\K\do\L\do\M\do\N%
	\do\O\do\P\do\Q\do\R\do\S\do\T\do\U\do\V\do\W\do\X%
	\do\Y\do\Z}

\begin{document}
%
\title{Optimal Combination of Frequency Control and Peak Shaving with Battery Storage Systems}

\author{Jonas~Engels,~\IEEEmembership{Student~Member,~IEEE,}
	Bert~Claessens 
	and~Geert~Deconinck,~\IEEEmembership{Senior~Member,~IEEE}
	
	\thanks{Jonas Engels is with Centrica Business Solutions Belgium NV, 2600 Antwerp, Belgium and with the Department of Electrical Engineering, KU Leuven/EnergyVille, Leuven, Belgium  (jonas.engels@restore.energy)}
	\thanks{Bert Claessens is with Centrica Business Solutions Belgium NV, 2600 Antwerp, Belgium (bert.claessens@restore.energy)}
	\thanks{Geert Deconinck is with the Department of Electrical Engineering, KU Leuven/EnergyVille,	Leuven, Belgium (geert.deconinck@kuleuven.be)}
	\thanks{This work is partially supported by Flanders Innovation \& Entrepreneurship (VLAIO) (grant no. BM 160214)}
	\thanks{\copyright 2019 IEEE. Personal use of this material is permitted.  Permission from IEEE must be obtained for all other uses, in any current or future media, including reprinting/republishing this material for advertising or promotional purposes, creating new collective works, for resale or redistribution to servers or lists, or reuse of any copyrighted component of this work in other works.}
\thanks{Published in IEEE Transactions on Smart Grid on 31 December 2019. DOI: 10.1109/TSG.2019.2963098}}


\maketitle

\begin{abstract}
Combining revenue streams by providing multiple services with battery storage systems increases profitability and enhances the investment case.
In this work, we present a novel optimisation and control framework that enables a storage system to optimally combine the provision of primary frequency control services with peak shaving of a consumption profile.

We adopt a dynamic programming framework to connect the daily bidding in frequency control markets with the longer term peak shaving objective: reducing the maximum consumption peak over an entire billing period. 
The framework also allows to aggregate frequency control capacity of multiple batteries installed at different sites, creating synergies when the consumption profile peaks occur on different times.

Using a case study of two batteries at two industrial sites, we show that the presented approach increases net profit of the batteries significantly compared to using the batteries for only peak shaving or frequency control.




\end{abstract}

\section{Introduction}\label{sec:intro}
\IEEEPARstart{B}{attery} energy storage systems (BESSs) installed behind-the-meter at the consumer's premises can be used for a variety of different services~\cite{MALHOTRA2016}. Often, the purpose of such a BESS is to decrease the energy costs of the consumer by optimising the charging schedule of the BESS towards their energy tariff. 
In case the consumer faces peak demand charges, usually a part of the grid tariffs, performing peak shaving with the BESS, i.e. reducing the consumer's power consumption peak, is an effective way to decrease energy costs~\cite{Wu2016}.

A BESS installed behind-the-meter can also be used to provide ancillary services, such as frequency control, to the transmission system operator (TSO). 
Especially primary frequency control (of frequency containment reserves) and frequency regulation services are seen to be a good match for a BESS, as the service provides a relatively high remuneration~\cite{Oudalov2006}, requires only a short duration of activation and a fast response, all of which a BESS can provide without problems~\cite{Zhang2017}.

By using the BESS for both energy tariff optimisation and frequency control service, one can combine both revenue streams, increase profitability and build a stronger business case for the investment in the BESS.
However, having a BESS providing both services concurrently is not straightforward from a control and optimisation perspective. 
One faces a trade-off, as using the BESS more for frequency control will decrease its peak shaving potential, which can be optimised. 

\subsection{Frequency Control with a BESS}\label{sec:fcr_bess}
The focus of this paper will be on primary frequency control services or frequency containment reserves (FCR), as defined by ENTSO-E~\cite{ENTSO-E2013}, since mainly the FCR markets have been opening up for third parties in Europe. 
Nevertheless, the presented methodology could also be applied to other types of frequency control or frequency regulation.

FCR is used by TSOs to stabilise the grid frequency continuously and thereby balance the demand and production in the grid.
When participating in FCR, one sells a certain amount of symmetric FCR capacity $r$ to the TSO which has to be available during the entire contracted period.
In FCR, one has to adjust its power $P^{FCR}$ proportionally to deviations of the grid frequency $f(t)$ from the nominal grid frequency $f_{nom}$ (\SI{50}{Hz} in Europe), so that the sold FCR capacity $r$ is reached when the frequency deviation is at a predefined maximum $\Delta f_{max}$ (= \SI{200}{mHz} in the Continental Europe synchronous region): $P^{FCR}(t) = r \Delta f(t) = r (f(t) - f_{nom})/\Delta f_{max}$.
In line with the recent changes in the German FCR market \textit{Regelleistung}~\cite{Regelleistung, ENTSOE_consult_FCR}, we assume a daily bidding process with daily auctions.

When delivering FCR with a BESS for a while, the battery can become full or empty at which point it is unable to provide the symmetric FCR service any further. Therefore, a state of charge (SoC) control strategy, or \emph{recharge} controller, which manages the SoC to ensure the BESS can deliver the FCR capacity for the entire contract period, is necessary~\cite{Fleer2016}.


\subsection{Peak Shaving}
Grid tariffs for commercial and industrial consumers usually consist of an energy charge (in \EUR/kWh) and a demand charge $c_{peak}$ (in \EUR/kW), where the latter is a charge proportional to the highest metered consumption peak during the billing period~\cite{Picciariello2015}. 
Such demand charges are typically used to recover the capacity-based costs of the network infrastructure, and are foreseen to become increasingly important with a growing share of distributed generation~\cite{Passey2017}.
With this tariff structure, a BESS can reduce network costs by discharging at the moments when the site is consuming its maximum power and charging when the site is consuming less, thereby reducing the site's metered consumption peak.

In practice, the highest metered consumption corresponds to the highest $n$-minute averaged power of the site, as usually energy meters with an $n$-minute resolution are used for settlement.
In this work, we consider demand charges proportional to the maximum quarter-hourly average power over one month, corresponding the German network tariff structure~\cite{BayernwerkNetzengelte}.






\subsection{Related Literature}
Previous work~\cite{LEADBETTER2012685, Koller2015, Bhattarai2016, OudalovPeakShaving} shows the ability of BESSs to perform peak shaving while making use of various control methods.
A simple threshold control policy is used in~\cite{LEADBETTER2012685}, while~\cite{Koller2015} demonstrates the use of a model predictive control to perform peak shaving of an office building.
A two-stage control methodology, consisting out of a scheduling stage and a real-time adaptive control stage
, is presented in~\cite{Bhattarai2016}, while \cite{OudalovPeakShaving} proposes dynamic programming to perform peak shaving.




Other work~\cite{Schiapparelli2018, Oudalov2007, Megel2013} has been devoted to BESSs providing frequency control services and the design of a recharge controller. 
The work in \cite{Oudalov2007} shows that a recharge controller is needed when using a BESS for FCR and propose a simple, threshold based recharge controller. 
More advanced, heuristic control strategies are proposed in~\cite{Megel2013}, which give better performance, but do not ensure any form of optimality. 

Few authors however have looked at the combination of both services.
Braeuer et al.~\cite{BRAEUER2019} perform a high-level economic analysis of BESSs combining peak shaving with frequency control. 
A similar approach is followed in~\cite{Moreno2015}, but with the peak shaving objective formulated as a hard network constraint, rather than implicitly through a demand charge.
Both papers indicate a significant added value in combining both services, however they assume perfect hindsight of the stochastic variables and do not develop a controller able to deliver the combination of services in day-to-day operation.


This work fills this gap by presenting 
an operational control framework
that enables a BESS 
to successfully combine peak shaving with frequency control services. 
The presented method extends our previous work on frequency control with BESSs~\cite{Engels2017}, by adding the peak shaving objective using dynamic programming and a customised stochastic optimisation objective. 
The main contributions of this paper are:
\begin{itemize}
	\item A novel stochastic optimisation and control framework that is able to optimally combine frequency control with peak shaving objectives using a BESS.
	\item A methodology which can be applied to efficiently
	aggregate frequency control capacity of multiple BESSs at different sites.
	\item A case study of two real industrial sites which shows that the presented approach increases value of the BESS compared to using the BESS for only a single objective.
\end{itemize}


In what follows, a bold symbol $\bm{x} = (x_0, \ldots, x_{n_t})^T$ denotes a vector containing the variables $x_k$, while $\E[\cdot]$ is the expected value operator, $\Pr(x\leq X)$ the probability of $x\leq X$, $\overline{X}$ the mean value of $X$, $[\cdot]^+~\equiv~\max(\cdot,0)$ operating element-wise on vectors, $\max(\bm{x})~\equiv~\max_k(x_k)$ the maximum element of the vector $\bm{x}$ and $I_{n_t} \in \mathbb{R}^{n_t \times n_t}$ the identity matrix.

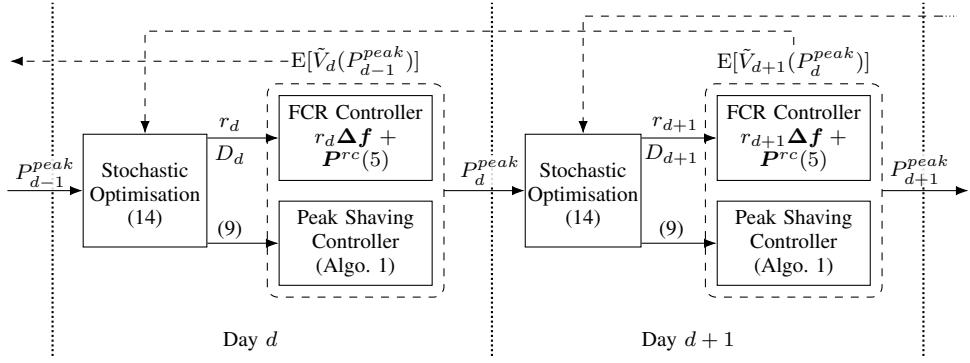
\begin{figure*}[t]
	\centering
	
	\begin{tikzpicture}
	\tikzstyle{every node}=[font= \footnotesize]
	
	\def \rect{0.15}
	\def \dScheduleControl{2.8}
	\def \widthControl{1.8cm}
	
	\node(schedule) [draw=black, align=center, minimum height=1.5cm, inner sep=2.5] at (-0.5,0)
	{\\Stochastic \\ Optimisation \\ (\ref{eq:DP})};
	
	\draw [thick, densely dotted] 
	($(schedule.west) + (-0.4,-2.2)$) -- ($(schedule.west) + (-0.4,2.5)$);

	\node [fill=white,text opacity=1,fill opacity=0.85, inner sep=0.5pt] at
	($(schedule.west) + (-0.5,0.23)$) {$P^{peak}_{d-1}$};
	
	\draw [-Latex] ($(schedule.west) + (-1.0,0)$) -- (schedule.west);
	
	\node(fcr) [draw=black, text width=\widthControl, text centered] at
	($(schedule) + (\dScheduleControl,0.7)$) 
	{FCR Controller \\ $r_d\bm{\Delta f} + \bm{P}^{rc}(\ref{eq:state_feedback})$ };
	\node(ps) [draw=black, text width=\widthControl, text centered] at 
	($(schedule) + (\dScheduleControl,-0.7)$) 
	{Peak Shaving Controller \\ (Algo. \ref{alg:rulebased})};
	
	\draw [-Latex] (schedule.east |- fcr) -- (fcr);
	\draw [-Latex] (schedule.east |- ps) -- (ps);
	
	\draw[rounded corners, dashed] ($(ps.south west) + (-\rect,-\rect)$) rectangle ($(fcr.north east) + (\rect,\rect)$);

	\node [fill=white,text opacity=1,fill opacity=0.85, inner sep=0.0pt] at ($(schedule.east |- fcr) + (0.3,0.18)$) {$r_d$};
	\node [fill=white,text opacity=1,fill opacity=0.85, inner sep=0.1pt] at ($(schedule.east |- fcr) + (0.3,-0.22)$) {$D_d$};
	\node [fill=white,text opacity=1,fill opacity=0.85, inner sep=0.0pt] at ($(schedule.east |- ps) + (0.3,0.18)$) {(\ref{eq:FCR_ps_div})};
	
	\node(scheduleD1) [draw=black, align=center, minimum height=1.5cm, inner sep=1] at	($(fcr.east)!0.5!(ps.east) + (2.,0)$) 
	{\\Stochastic \\ Optimisation \\ (\ref{eq:DP})};

	\node(fcrD1) [draw=black, text width=\widthControl, text centered] at
	($(scheduleD1) + (\dScheduleControl,0.7)$) 
	{FCR Controller \\ $r_{d+1}\bm{\Delta f} + \bm{P}^{rc}(\ref{eq:state_feedback})$ };
	\node(psD1) [draw=black, text width=\widthControl, text centered] at 
	($(scheduleD1) + (\dScheduleControl,-0.7)$) 
	{Peak Shaving Controller \\ (Algo. \ref{alg:rulebased})};
	
	\draw [-Latex] (scheduleD1.east |- fcrD1) -- (fcrD1);
	\draw [-Latex] (scheduleD1.east |- psD1) -- (psD1);
	\draw[rounded corners, dashed] 
	($(psD1.south west) + (-\rect,-\rect)$) rectangle ($(fcrD1.north east) + (\rect,\rect)$);
	
	\node [fill=white,text opacity=1,fill opacity=0.85, inner sep=0.0pt] at ($(scheduleD1.east |- fcrD1) + (0.45,0.18)$) {$r_{d+1}$};
	\node [fill=white,text opacity=1,fill opacity=0.85, inner sep=0.0pt] at ($(scheduleD1.east |- fcrD1) + (0.40,-0.22)$) {$D_{d+1}$};
	\node [fill=white,text opacity=1,fill opacity=0.85, inner sep=0.0pt] at ($(scheduleD1.east |- psD1) + (0.40,0.18)$) {(\ref{eq:FCR_ps_div})};
	
	\draw [-Latex] ($(fcr.east)!0.5!(ps.east) + (\rect,0)$) -- (scheduleD1.west);
	
	\draw [thick, densely dotted] ($(fcr)!0.5!(ps)!0.5!(scheduleD1) + (0.3,-2.2)$) -- ($(fcr)!0.5!(ps)!0.5!(scheduleD1) + (0.3,2.5)$);
	
	\node [fill=white,text opacity=1,fill opacity=0.85, inner sep=0.10pt] at 
	($(fcr.east)!0.5!(ps.east)!0.5!(scheduleD1.west)+ (0.15,0.27)$)
	{$P^{peak}_{d}$};
	
	\node(V) [fill=white,text opacity=1,fill opacity=0.85, inner sep=0.75pt] at
	($(fcr.north)+(0,0.45)$) 
	{$\E[\tilde{V}_{d}(P^{peak}_{d-1})]$};
	
	\draw [dashed,-Latex] (V) -- ($(V.west) + (-3.7,0)$);
	
	\draw [dashed,-Latex] ($(V.east) + (6.9,0.6)$)  -| (scheduleD1);
	\draw [densely dotted] ($(V.east) + (6.9,0.6)$) -- ($(V.east) + (7.1,0.6)$);
	
	\node(Vnew) [fill=white,text opacity=1,fill opacity=0.9, inner sep=0.5pt] at ($(fcrD1.north)+(0,0.45)$) 
	{$\E[\tilde{V}_{d+1}(P^{peak}_{d})]$};
	
	\draw [dashed,-Latex] (Vnew) -- +(0,0.4) -| (schedule);
	
	\path let \p1 = ($(fcrD1)!0.5!(scheduleD1)$) in node at (\x1,-2.0) 
	{Day $d+1$};
	\path let \p1 = ($(fcr)!0.5!(schedule)$) in node  at (\x1,-2.0) 
	{Day $d$};
	
	\draw [thick, densely dotted] 
	($(fcrD1.east)!0.5!(psD1.east) + (0.7,-2.2)$) -- ($(fcrD1.east)!0.5!(psD1.east) + (0.7,2.5)$);
	
	\draw [-Latex] 
	($(fcrD1.east)!0.5!(psD1.east) + (\rect,0) $) -- ($(fcrD1.east)!0.5!(psD1.east) + (1.3,0)$);
	
	\node [fill=white,text opacity=1,fill opacity=0.85, inner sep=0.5pt] at
	($(fcrD1.east)!0.5!(psD1.east) + (\rect+0.5,0.27) $) 
	{$P^{peak}_{d+1}$};

	\end{tikzpicture}
	\caption{Overview of the dynamic programming-based optimisation and control framework to combine peak shaving and frequency control with a BESS.}
	\label{fig:overview}
\end{figure*}

\section{Optimisation and Control Framework}\label{sec:prob}

A schematic overview of the optimisation and control framework to combine peak shaving and frequency control with a BESS at a particular site is shown in Figure \ref{fig:overview}.
In an FCR market with daily auctions, one has to decide every day $d$ on the FCR capacity $r_d$ one wants to sell. In the proposed framework, we make this decision through a stochastic optimisation problem (\ref{eq:DP}). The results of this optimisation are then used in the real-time FCR (\ref{eq:state_feedback}) and peak shaving (Algo. \ref{alg:rulebased}) controllers that are able to deliver both applications in real-time. 
As the peak shaving objective concerns the maximum consumption peak over an entire month, dynamic programming is adopted using the observed peak power $P^{peak}_d$ at the end of day $d$ as the state in the stochastic optimisation of the next day $d+1$.

Section~\ref{sec:freqcontrol} first describes the optimisation and control of a BESS for frequency control during one day. Subsequently, Section~\ref{sec:combining} explains how we add the peak shaving objective in the optimisation problem and proposes a real-time peak shaving controller, which can be used in practice to perform peak shaving in real-time, i.e. without the need for any hindsight knowledge. 
The section also elaborates how we extend the peak shaving objective from the maximum consumption peak during one day towards one month using dynamic programming and value function approximation.



In the control framework, we employ a BESS model with constant charging and discharging efficiencies $\eta^c, \eta^d$, discretised with time step $\Delta t$:
\begin{IEEEeqnarray}{rCCCl}\IEEEyesnumber\label{eq:bat_model}
E_{min} &\leq& E_k^{bat}& \leq& E_{max}, \quad P_{min} \leq P_k^{bat} \leq P_{max},\IEEEyessubnumber\\
E_{k+1}^{bat} &=& \IEEEeqnarraymulticol{3}{l}{ E_{k}^{bat} + (\eta^{c} [P_{k}^{bat}]^+ - \frac{1}{\eta^{d}}  [-P_{k}^{bat}]^+) \Delta t,} \IEEEyessubnumber \label{eq:bat_model_energy}
\end{IEEEeqnarray}
where $P_k^{bat}, E_k^{bat}$ are the power and energy content of the BESS at time step $k$ respectively.

%


\section{Frequency Control Framework}\label{sec:freqcontrol}
The frequency control framework we use in this paper is an extension of the robust optimisation presented in our previous work~\cite{Engels2017}, which we will shortly summarise here. For detailed information, the reader is referred to the original paper.

The goal here is to determine both the maximum amount of frequency control capacity $r$ the BESS can provide during one day and the power of the recharge controller $\bm{P}^{rc} \in \mathbb{R}^{n_t}$ that ensures the BESS stays within its SoC boundaries when delivering the FCR service. 
The electricity costs for charging and discharging the battery will not be taken into account here, but will be added in Section~\ref{sec:combining}. 
As the frequency deviation profile $\bm{\Delta f} \in \mathbb{R}^{n_t}$ is inherently stochastic, the energy $\bm{E}^{bat}$ and recharge power $\bm{P}^{rc}(\bm{\Delta f})$, which are dependent on the frequency profile, are also stochastic. 
The optimisation, maximising revenues from providing frequency control capacity $r$ at a price $c^{FCR}$, can then be formulated as a chance-constrained problem:
\begin{mini!}[0]
{}{-c^{FCR} r } 
{\label{eq:problem_r1}}{}
\addConstraint{\bm{P}^{bat}}{= r\bm{\Delta f} + \bm{P}^{rc}}
\addConstraint{1-\epsilon}{\leq \Pr( E_{min} \leq \bm{E}^{bat}) \label{constr:Emin_prob}}
\addConstraint{1-\epsilon}{\leq \Pr( \bm{E}^{bat} \leq E_{max}) \label{constr:Emax_prob}}
\addConstraint{1-\epsilon}{\leq \Pr(P_{min}+r \leq \bm{P}^{rc}) \label{constr:Pmin_prob}}
\addConstraint{ 1-\epsilon}{\leq \Pr(\bm{P}^{rc} \leq {P}_{max}-r) \label{constr:Pmax_prob}}
\addConstraint{	E_{k+1}^{bat}}{ = E_{k}^{bat} + (\eta^{c} [P_{k}^{bat}]^+ - \frac{1}{\eta^{d}}  [-P_{k}^{bat}]^+) \Delta t. \label{constr:energy_r1}}
\end{mini!}
The chance constraints (\ref{constr:Emin_prob})-(\ref{constr:Pmax_prob}) arise here naturally as the energy $\bm{E}^{bat}$ and recharge power $\bm{P}^{rc}$, which are stochastic due to their dependency on the frequency, have to be limited to the energy and power capacity of the BESS.

We solve (\ref{eq:problem_r1}) using robust optimisation~\cite{Ben-Tal2009}, as it generates a safe approximation to (\ref{constr:Emin_prob})-(\ref{constr:Pmax_prob}) while allowing to make $\epsilon$ arbitrary small and retaining a tractable and efficiently solvable second-order cone problem (SOCP). To achieve this, a couple of reformulations are needed. 


\subsection{Battery Efficiency}
To avoid the integer variables resulting from the $[\cdot]^+$ operators in (\ref{constr:energy_r1}), we set the efficiencies in the constraint (\ref{constr:energy_r1}) itself to $\eta^c=\eta^d=1$. In turn, we incorporate the efficiencies into the frequency deviations when discretising them:
\begin{equation}\label{eq:en_freq}
\Delta f_k = \frac{1}{\Delta t}\int_{(k-1)\Delta t}^{k\Delta t} \Big(\eta^c\big[\Delta f(t)\big]^+ - \frac{1}{\eta^d}\big[-\Delta f(t)\big]^+\Big) dt.
\end{equation} 
In our previous work~\cite{Engels2017}, we showed that this approximation does not lead to violations of the constraints when $\eta^c,\eta^d <1$. 


\subsection{Recharge Controller}
The recharge power $\bm{P}_{rc}$ in (\ref{eq:problem_r1}) has to ensure the probabilities of (\ref{constr:Emin_prob})-(\ref{constr:Pmax_prob}) are satisfied. As the frequency deviations $\bm{\Delta f}$ are gradually revealed over time, 
we can have $\bm{P}_{rc}$ be dependent on the $n_{rc}$ past frequency deviations: $P^{rc}_k = \pi^k(\Delta f_{\left[k-n_{rc}\right]^+}, \ldots, \Delta f_{k-1})$, with $\pi^k$ a policy at time step $k$.
As an optimisation over functions $\pi^k$ is generally intractable, we limit ourselves to a linear policy:
\begin{equation}\label{eq:recharge_policy}
P_k^{rc} = \sum_{i=\left[k-n_{rc}\right]^+}^{k-1} d_{ki}\Delta f_i, \quad \bm{P}^{rc} = D\bm{\Delta f},  
\end{equation}
with $d_{ki}$ the coefficients of the recharge strategy, contained in the lower triangular matrix $D \in \mathbb{R}^{n_t \times n_t}$. 

As this recharging policy will be calculated with the efficiencies incorporated in the frequency signal (\ref{eq:en_freq}) and not in the battery model itself, the policy will not be directly applicable to a real battery system with $\eta^c,\eta^d <1$. However, following~\cite{Goulart2006b}, such a linear disturbance feedback policy can be transformed into an equivalent state-feedback policy:
\begin{equation}\label{eq:state_feedback}
\bm{P^{rc}} = (I+\frac{1}{r}D)^{-1} \frac{1}{r} D\bm{\Delta E^{bat}},
\end{equation}
with $\Delta E^{bat}_k = (E_k^{bat}-E_{k-1}^{bat})/\Delta t$. In this form, the recharge power reacts on the past states, which include the effect of the actual efficiency losses and other non-linearities of the BESS.
The FCR controller of the BESS is then: $\bm{P}^{bat} = r\bm{\Delta f} + \bm{P}^{rc}$.

\subsection{Robust Reformulation}
With the adaptations described above, we can use robust optimisation to create a safe approximation of the chance constraints (\ref{constr:Emin_prob})-(\ref{constr:Pmax_prob}). 
The idea is to design an uncertainty set of frequency deviations $\bm{\Delta f} \in \mathcal{F}_\epsilon$, against which each of the constraints (\ref{constr:Emin_prob})-(\ref{constr:Pmax_prob}) have to be satisfied at all times:
\begin{equation}\label{eq:RO}
\max_{\bm{\Delta f} \in \mathcal{F}_\epsilon} \bm{a_i}^T \bm{\Delta f} \leq b_i, \qquad i = 1,\ldots,n_c,
\end{equation}
with $(\bm{a_i}, b_i)$ defined to represent one inequality in the probability operators of (\ref{constr:Emin_prob})-(\ref{constr:Pmax_prob}) and $n_c = 4 n_t$. 
Let $A = [-(D+rI)^T G^T | (D+rI)^T  G^T | -D^T | D^T ]^T$, with $G$ a lower triangular matrix with $\Delta t$ as elements, and vector $\bm{b} = [-(E_{min}-E^{bat}_0) | E_{max}-E^{bat}_0 | -(P_{min}+r) | P_{max}-r ] ^T$, then $\bm{a_i}^T$ is the $i$-th row of $A$ and $b_i$ the $i$-th element of $\bm{b}$.

Chen et al. show in \cite{Chen2009} that an asymmetric uncertainty set based on forward $\sigma_{f_k}(\Delta f_k)$ and backward $\sigma_{b_k}(\Delta f_k)$ deviations, which can be estimated from samples of $\Delta f_k$, provides the tightest bound for small $\epsilon$.
With $Q = \text{diag}(\sigma_{f1},\ldots,\sigma_{fn_t} )$ and $R = \text{diag}(\sigma_{b1},\ldots,\sigma_{bn_t})$, we can reformulate the constraints~(\ref{eq:RO}) into the following second-order cone constraints:
\begin{equation}\label{eq:R1_constr}
\bm{a_i}^T \overline{\bm{\Delta f}} + \sqrt{-2\ln{\epsilon}} 
\lVert \bm{u_i}  \rVert_2
\leq b_i, \qquad i = 1,\ldots,n_c ,
\end{equation}
where $\bm{u_i} = \max(Q\bm{a_i}^T W^{-1},-R\bm{a_i}^T W^{-1})$, with the maximum taken element-wise and $W^T W = \Sigma_{\bm{\Delta f}}^{-1}$ the Cholensky decomposition of the inverse of the covariance matrix $\Sigma_{\bm{\Delta f}}$ of ${\bm{\Delta f}}$. 
We refer to our previous work~\cite{Engels2017} for the details on the derivation of (\ref{eq:R1_constr}).

With these reformulations, (\ref{eq:problem_r1}) becomes the following tractable second-order cone problem, which can be readily solved by various commercial and non-commercial solvers:
\begin{mini!}[0]
{}{-c^{FCR} r } 
{\label{eq:problem_r1_ref}}{}
\addConstraint{\bm{a_i}^T \overline{\bm{\Delta f}} + \sqrt{-2\ln{\epsilon}} 
	\lVert \bm{u_i} \rVert_2}{\leq b_i, \label{constr:problem_r1_ref_soc}}{\quad i = 1,\ldots,n_c}
\addConstraint{Q\bm{a_i}^T W^{-1}}{\leq \bm{u_i},  \label{constr:problem_r1_ref_uq}}{\quad i = 1,\ldots,n_c}
\addConstraint{-R\bm{a_i}^T W^{-1}}{\leq \bm{u_i}, \label{constr:problem_r1_ref_ur}}{\quad i = 1,\ldots,n_c.}
\end{mini!}

\section{Combining Peak Shaving and Frequency Control}\label{sec:combining}
When adding the peak shaving objective to the optimisation~(\ref{eq:problem_r1}), one has to ensure the chance constraints (\ref{constr:Emin_prob})-(\ref{constr:Pmax_prob}) are still satisfied.
To achieve this, we split the BESS into two virtual batteries: one for peak shaving and one for frequency control. By constraining the virtual battery for frequency control to (\ref{eq:R1_constr}), it is ensured (\ref{constr:Emin_prob})-(\ref{constr:Pmax_prob}) are satisfied.
Besides, by intelligently shaping the virtual battery boundaries, one can obtain synergies. For instance, one can reserve less recharge power and hence more power for peak shaving at the moments when consumption peaks are expected, and compensate for this at the moments where consumption is expected to be low.

For a specific FCR capacity $r$ and recharge policy $D$, equation (\ref{eq:RO}) allows to calculate the minimum and maximum power $(\bm{P}_{min}^{FCR}, \bm{P}_{max}^{FCR})$ and energy $(\bm{E}_{min}^{FCR}, \bm{E}_{max}^{FCR})$ capacity needed to perform frequency control at any time step $k$. 
The remaining power and energy capacity of the BESS can then be used to perform peak shaving: 
\begin{IEEEeqnarray}{rCl?rCl}
\bm{P}_{min}^{ps} &=& P_{min} - \bm{P}_{min}^{FCR},& \bm{P}_{max}^{ps} &=& P_{max} - \bm{P}_{max}^{FCR}, \IEEEyesnumber \label{eq:FCR_ps_div} \IEEEeqnarraynumspace \IEEEyessubnumber\\
\bm{E}_{min}^{ps} &=& E_{min} - \bm{E}_{min}^{FCR},& \bm{E}_{max}^{ps} &=& E_{max} - \bm{E}_{max}^{FCR}, \IEEEyessubnumber
\end{IEEEeqnarray}


Let $\bm{P}^{ps}$ and $\bm{E}^{ps}$ be the power and energy profile of the part of the BESS used for peak shaving and $\bm{P}^{prof}$ the consumption profile of the site. The combined optimisation, maximising frequency control revenues, minimising the expected peak power costs of the site and the additional electricity cost due to the charging of the battery, can then be formulated as:
\begin{IEEEeqnarray}{t?rCl}
min& \IEEEeqnarraymulticol{3}{l}{\E[{c}^{peak} P^{peak} + C^{elec} ] - c^{FCR} r,} \IEEEyesnumber \label{eq:problem} \IEEEyessubnumber \label{opt:prob_obj} \\
\textrm{s.t.}&P^{peak} & = & \max(P^{grid}_0, \ldots, P^{grid}_{n_t}), \IEEEyessubnumber\\
&\bm{P}^{grid} & = &\bm{P}^{prof} + \bm{P}^{ps} + (D + rI_{n_t}) \bm{\Delta f}, \IEEEyessubnumber\\
&C^{elec} &=&c_{elec}  \left(\bm{P}^{ps} + \left(D + r I_{n_t}\right) \bm{\Delta f} \right)\Delta t, \IEEEyessubnumber \\
&P_{min} &\leq& \bm{P}_{min}^{ps} \leq \bm{P}^{ps} \leq \bm{P}_{max}^{ps} \leq P_{max}, \IEEEyessubnumber \label{eq:ps_p_constraints}  \\
&E_{min} &\leq& \bm{E}_{min}^{ps} \leq \bm{E}^{ps} \leq \bm{E}_{max}^{ps} \leq E_{max}, \IEEEyessubnumber \label{eq:ps_e_constraints}\\
&\IEEEeqnarraymulticol{3}{l}{(\text{\ref{eq:bat_model_energy}), (\ref{constr:problem_r1_ref_soc}),  (\ref{constr:problem_r1_ref_uq}), (\ref{constr:problem_r1_ref_ur}), (\ref{eq:FCR_ps_div})}, } \IEEEyessubnumber \label{eq:prob_r1_constr}
\end{IEEEeqnarray} 
with $c_{elec}$ the per unit energy cost.

\subsection{Stochastic Optimisation}\label{sec:stoch_opt}
The expected value operator in the objective (\ref{opt:prob_obj}) depends on the stochastic consumption profile $\bm{P}^{prof}$ and frequency deviation  profile $\bm{\Delta f}$ and thus concerns an $n_t$-dimensional integration, which is intractable in practice.
To approximate the expected value operator, one can use a Sample Average Approximation (SAA)~\cite{Shapiro2009} by taking the empirical mean over independent and identically distributed (iid) samples or \emph{scenarios} of the stochastic variables. 
With $\bm{P}^{prof}_j$,  $\bm{\Delta f}_j$ the $j$-th iid consumption profile and frequency deviation sample respectively, $j= 1,\dots, n_{sc}$, and $p_j= 1/n_{sc}$ the probability of scenario $j$, one can approximate the expected value operator as follows:
\begin{IEEEeqnarray}{lCl}
\IEEEeqnarraymulticol{3}{l}{\E[c_{peak} P^{peak} + C^{elec}] \approx \sum_{j=1}^{n_{sc}} p_j \left({c}_{peak} P^{peak}_j + C^{elec}_j \right),}  \IEEEyesnumber\label{eq:stoch_opt} \IEEEeqnarraynumspace \IEEEyessubnumber \label{eq:stoch_opt_obj}\vspace{-0.3cm}
\end{IEEEeqnarray}
where:
\begin{IEEEeqnarray}{lCl}
P^{peak}_j &=& \max \left(\bm{P}^{prof}_j + \bm{P}^{ps}_j + \left(D + r I_{n_t}\right) \bm{\Delta f}_j \right),\IEEEyessubnumber \IEEEeqnarraynumspace \label{eq:peak_stochopt}\\
C^{elec}_j &=&c_{elec}  \left(\bm{P}^{ps}_j + \left(D + r I_{n_t}\right) \bm{\Delta f}_j \right)\Delta t. \IEEEyessubnumber \label{eq:c_elec}
\end{IEEEeqnarray}

\subsubsection{Interference Peak Shaving and Frequency Control}
In case a positive frequency control power is required ($\Delta f_{k,j}>0$) when the consumption of the site is high, this could increase the peak consumption $P^{peak}$ of the site. 
Using (\ref{eq:stoch_opt}) in the optimisation problem (\ref{eq:problem}) would then result in a peak shaving power $P^{ps}_{k,j}$ which completely compensates for the frequency control power: $P^{ps}_{k,j} = -r \Delta f_{k,j}$. This means that in practice, no frequency control power has been delivered to the grid.

To prevent this interference in the optimisation, the peak shaving power $\bm{P}^{ps}_j$ should be independent of the required frequency control power $r \bm{\Delta f}_{j}$. 
We achieve this by sampling the frequency profile $\bm{\Delta f}$ separately from the consumption profile $\bm{P}^{prof}$
and have each peak shaving power scenario $\bm{P}^{ps}_j$ dealing with all frequency deviation samples. Let $v = 1,\ldots,n_v$ and $w = 1,\ldots,n_w$ be the index of the consumption profile samples $\bm{P}^{prof}_v$ and frequency deviation samples $\bm{\Delta f}_{w}$, respectively, then:
\begin{IEEEeqnarray}{rCl}\IEEEyesnumber
P^{peak}_{j} &=& \max \left(\bm{P}^{prof}_{v} + \bm{P}^{ps}_{v} + \left(D + r I_{n_t}\right) \bm{\Delta f}_{w} \right), \label{eq:combfreq_cons} \IEEEyessubnumber \\
C^{elec}_{j} &=&c_{elec}  \left(\bm{P}^{ps}_{v} + \left(D + r I_{n_t}\right) \bm{\Delta f}_{w} \right),  \IEEEyessubnumber\\
p_{j} & =& p_{v}^{prof} p_{w}^{\Delta f}, \IEEEyessubnumber\\
j &\coloneqq& v n_{w} + w, \quad v = 1,\ldots,n_{v}, \quad w= 1,\ldots,n_{w}. \nonumber 
\end{IEEEeqnarray}
Here, $j = 1,\ldots n_{sc}$, with $n_{sc}=n_{w} n_{v}$, is the index used in the SAA objective (\ref{eq:stoch_opt_obj}), $p_{v}^{prof}$ the probability of the consumption profile scenario $w$ and $p_{w}^{\Delta f}$ the probability of the frequency deviation scenario $v$.
With this approach, each battery peak shaving power scenario $\bm{P}^{ps}_{v}$ is able to reduce the peak of the corresponding consumption profile $\bm{P}^{prof}_{v}$, but also has to deal with 
all frequency deviation profiles $\bm{\Delta f}_{w}$ in the optimisation.

\subsubsection{Scenario Reduction}
As the SAA objective (\ref{eq:stoch_opt_obj}) converges to the true value with a rate of $O(1/{n_{sc}})$~\cite{Shapiro2009}, a high number of scenarios are needed to reach an acceptable accuracy.
To reduce the number of scenarios and increase computational efficiency, we employ the \emph{fast forward selection} algorithm presented by Heitsch and R\"{o}misch~\cite{Heitsch2003}.
The original fast forwards selection algorithm is a heuristic to minimise the Kantorovich distance $D_K(\Omega, \Omega_r)$ between an original set of scenarios $\Omega$ and a new, reduced set of scenarios $\Omega_r \subset \Omega$:
\begin{equation}\label{eq:scenred}
D_K(\Omega, \Omega_r) = \sum_{\bm{\omega} \in \Omega \setminus  \Omega_r} p_{\bm{\omega}} \min_{\bm{\omega}' \in \Omega_r} c(\bm{\omega}, \bm{\omega}'),
\end{equation}
with $p_{\bm{\omega}}$ the probability of scenario $\bm{\omega}$ and the cost function $c(\bm{\omega}, \bm{\omega}') = \left\|\bm{\omega} - \bm{\omega}' \right\|_2$~\cite{Romisch2002Stability}.


%
%
%

\begin{figure}
	\centering
	\subfloat[{$\E[P^{peak}]$} of the reduced set of consumption profile samples. \label{fig:scen_red_peak}]{%
		\includegraphics[width=0.475\columnwidth]{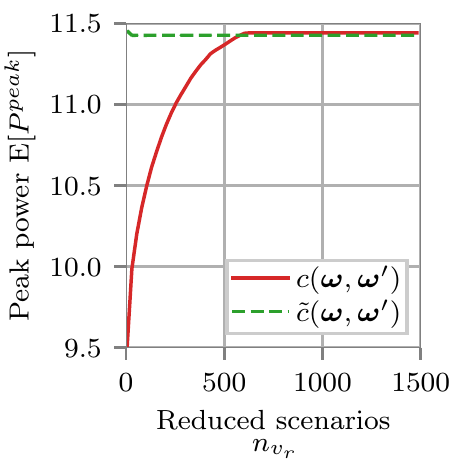}
	}\hfill
	\subfloat[SAA optimality gap when solving problem (\ref{eq:problem}) with $n_{v_r}=n_{w_r}$. \label{fig:scen_red_SAAgap}]{%
		\centering
		\includegraphics[width=0.475\columnwidth]{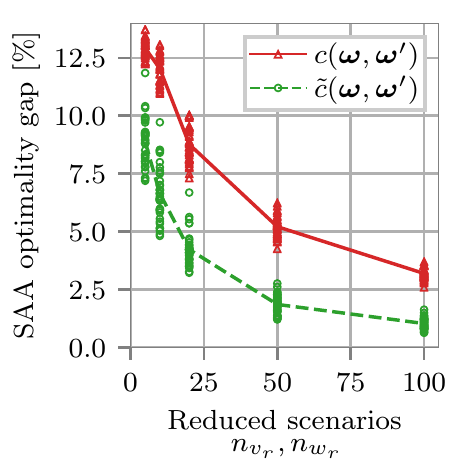}
	}
	\caption{(a) When reducing 1500 idd consumption profile samples to $n_{v_r}$ scenarios, the generic cost function $c(\bm{\omega}, \bm{\omega}')$ introduces a significant bias in $\E[P^{peak}]$ which the proposed cost function $\tilde{c}(\bm{\omega}, \bm{\omega}')$ is able to eliminate. (b)~The proposed cost function $\tilde{c}(\bm{\omega}, \bm{\omega}')$ also reduces the SAA optimality gap when solving problem (\ref{eq:problem}) with $n_{v_r}=n_{w_r}$ reduced scenarios.}
	\label{fig:scen_red}
\end{figure}

Figure \ref{fig:scen_red_peak} shows the expected peak power $\E[P_{peak}]$ during one day when reducing the number of scenarios 
using the generic cost function $c(\bm{\omega}, \bm{\omega}')$ in the fast forward selection algorithm. 
As one can see, the method introduces a significant bias when reducing to less than 600 scenarios. 
It has been noted previously in the literature~\cite{Morales2009} 
that using a cost function that is better able to capture the effect of adding a scenario on the objective of the problem can improve performance.
As our objective function (\ref{eq:stoch_opt_obj}) involves the maximum value of a scenario, we propose the following cost function in (\ref{eq:scenred}):
\begin{equation*}
\tilde{c}(\bm{\omega}, \bm{\omega}') = \left| \max(\bm{\omega}) - \max(\bm{\omega}')\right|.
\end{equation*}
As the dashed line in Figure~\ref{fig:scen_red_peak} shows, this new cost function is able to eliminate the bias on $\E[P_{peak}]$ 
almost completely.

To prevent the interference between peak shaving and frequency control objectives discussed above, we separately sample the consumption profiles $\bm{P}^{prof}_v$ and frequency deviations $\bm{\Delta f}_w$ and reduce them to $\bm{P}^{prof}_{v_r}$ with probability $p_{v_r}$, $v_r = 1,\ldots, n_{v_r}$  and $\bm{\Delta f}_{w_r}$ with probability $p_{w_r}$, $w_r = 1,\ldots, n_{w_r}$, respectively. 
We then combine the reduced scenarios as in (\ref{eq:combfreq_cons}), so that in total $n_{sc_r} = n_{v_r} n_{w_r}$ and $p_{j_r}=p_{v_r} p_{w_r}$ in (\ref{eq:stoch_opt_obj}).



Finally, Figure \ref{fig:scen_red_SAAgap} shows the optimality gap of  (\ref{eq:stoch_opt}) due to the SAA, calculated according to~\cite{Mak1999}, when reducing 1500 iid consumption and frequency deviation samples to $n_{v_r}$ and $n_{w_r}$ scenarios. 
The proposed cost function $\tilde{c}(\bm{\omega}, \bm{\omega}')$ decreases the SAA optimality gap with around \SI{50}{\%} for a same number of reduced scenarios, increasing computational efficiency.

\subsection{Non-Anticipative Peak Shaving Controller}\label{sec:pscontrol}
When solving the stochastic optimisation problem (\ref{eq:problem}) as described in the paragraphs above, one actually solves a two-stage stochastic problem. 
In a first stage, one decides on the FCR capacity $r$, the recharge policy $D$ and the peak shaving boundaries $(\bm{E}_{min}^{ps}, \bm{E}_{max}^{ps}, \bm{P}_{min}^{ps}, \bm{P}_{max}^{ps})$. 
In the second stage, one optimises the peak shaving power $\bm{P}^{ps}_{v_r}$ with complete (perfect hindsight) knowledge of the consumption profile $\bm{P}^{prof}_{v_r}$. 
In reality however, the consumption profile is only gradually revealed over time. Hence, the second stage solution is not usable in practice and a non-anticipative, potentially suboptimal, peak shaving control algorithm will be required.



Examples of such controllers vary from simple, rule-based controllers~\cite{LEADBETTER2012685} to model predictive control~\cite{Koller2015} and more complex dynamic programming methods~\cite{OudalovPeakShaving}.
The optimisation and control framework proposed in this paper allows the use of any of these control algorithms. 
However, to limit the scope of this paper we restrict ourselves to a rather simple, parametrised rule-based peak shaving policy.



Algorithm \ref{alg:rulebased} shows the proposed rule-based peak shaving controller. The controller discharges the battery every time $k$ the grid power $\hat{P}^{grid}_k$ surpasses a threshold $P_{thr}$ and recharges the battery every time $\hat{P}^{grid}_k$ goes below this threshold. 

\begin{algorithm}
	\caption{Rule-based peak shaving controller}
	\label{alg:rulebased}
	\begin{algorithmic}[1]
		\Parameter $P_{thr}^{init}, z_{\sigma}$.
		\Require $r, D, \bm{E}_{max}^{ps}, \bm{E}_{min}^{ps}, \bm{P}_{max}^{ps}, \bm{P}_{min}^{ps}$.
		\State $P_{thr} \gets P_{thr}^{init}$.
		\State $\overline{\bm{P}^{FCR}} \gets (D + r I) \overline{\bm{\Delta f}}$. 
		\State $\bm{s}^{P^{FCR}} \gets std\big[(D + r I) \bm{\Delta f}\big]$. 
		\ForAll{time step $k = 1 \ldots n_t$} \label{alg_s:iter_k}
		\State $\hat{P}^{grid}_k \gets P^{prof}_k + \overline{P^{FCR}_k} + z_\sigma s^{P_{FCR}}_k$. \label{alg_s:prof_add_fcr}
		\State $P^{ps}_k \gets P_{thr} - \hat{P}^{grid}_k$. \label{alg_s:add_thresh}
		\If{$P^{ps}_k < 0$} \Comment{Discharge Battery}
		\State $P^{ps}_k \gets \max(P^{ps}_k , P^{ps}_{min,k}, \eta_d (E^{ps}_{min,{k+1}} - E^{ps}_k)\Delta t )$,\label{alg_s:compare_discharge}
		\Else   \Comment{Charge Battery}
		\State $P^{ps}_k \gets \min(P^{ps}_k , P^{ps}_{max,k}, \frac{1}{\eta_c} (E^{ps}_{max,{k+1}} - E^{ps}_k)\Delta t )$.\label{alg_s:compare_charge}
		\EndIf
		\State $E^{ps}_{k+1} \gets E^{ps}_{k+1} + \eta_c \left[P^{ps}_k \right]^+ - \frac{1}{\eta_d}\left[-P^{ps}_k \right]^+$.
		\If{$\hat{P}^{grid}_k + P^{ps}_k > P_{thr} $} \Comment{Threshold Surpassed}
		\State $P_{thr} \gets \hat{P}^{grid}_k + P^{ps}_k$. \label{alg_s:update_thresh}
		\EndIf
		\EndFor
	\end{algorithmic}
\end{algorithm}

The battery power due to the frequency control $(D + r I) \bm{\Delta f}$ can induce additional power peaks, which we want to avoid as much as possible without hampering the actual FCR delivery. 
Therefore, in step \ref{alg_s:prof_add_fcr}, we compute a statistic of the FCR power to be delivered: the average FCR power $\overline{P^{FCR}_k}$ plus a factor $z_\sigma$ times the standard deviation of the FCR power $s^{P_{FCR}}_k$, which we add to the consumption profile $P^{prof}_k$ to obtain a modified grid power profile $\hat{P}^{grid}_k$, which is compared with the threshold $P_{thr}$ in step \ref{alg_s:add_thresh}.

Steps \ref{alg_s:compare_discharge} and \ref{alg_s:compare_charge} ensure that $\bm{P}^{ps}, \bm{E}^{ps}$ stay within the peak shaving boundaries $(\bm{E}_{min}^{ps}, \bm{E}_{max}^{ps}), (\bm{P}_{min}^{ps}, \bm{P}_{max}^{ps})$.
Finally, step \ref{alg_s:update_thresh} updates the threshold $P_{thr}$ if the battery was unable to keep the modified grid power $\hat{P}^{grid}_k $ below the threshold $P_{thr}$.

Algorithm \ref{alg:rulebased} has two parameters that can be freely chosen: $P_{thr}^{init}$ and $z_{\sigma}$, which and be used to adapt the controller to a specific configuration. 
For a particular value of these parameters, the performance of the controller can be evaluated by simulating the controller for a large number of iid consumption and frequency samples $n_{eval} \gg n_{sc}$, calculating the objective (\ref{eq:stoch_opt_obj}) and taking the empirical average over all scenarios.
To find the optimum values ${P_{thr}^{init}}^*, z_{\sigma}^*$, we then use a simple grid search.


\subsection{Dynamic Programming Framework}



The optimisation~(\ref{eq:problem}) considered so far deals with the daily decision making required in the FCR market. 
However, peak demand charges look at the highest peak over an entire billing period, here one month. To deal with these different time scales, we adopt a dynamic programming framework. Starting at the end of the month, we calculate the value of the objective $V_d(P^{peak}_{d-1})$ for each day $d= 1, \ldots, n_d$ of the month in function of the peak power $P^{peak}_{d-1}$ observed until the end of the previous day $d-1$. 
The daily optimisation becomes then:
\begin{IEEEeqnarray}{ll,rCl}
\IEEEeqnarraymulticol{5}{l}{V_d(P^{peak}_{d-1})}\nonumber \\* \quad 
= \> &\min& \IEEEeqnarraymulticol{3}{l}{\E[V_{d+1}(P^{peak}_d) + {C}^{elec}_d]-c^{FCR}_d r_d, } \IEEEeqnarraynumspace \, \IEEEyesnumber\label{eq:DP}  \IEEEyessubnumber \label{eq:DP_obj} \\
&\mathrm{s.t.}&P^{peak}_d &=& \max(P_{0,d}^{grid}, \ldots, P_{n_t,d}^{grid}, P^{peak}_{d-1}),  \IEEEyessubnumber \\
&&\bm{P}^{grid}_d &=& \bm{P}^{prof}_d + \bm{P}^{ps}_d + \left(D_d + r_d I_{n_t}\right) \bm{\Delta f}_d, \IEEEeqnarraynumspace \IEEEyessubnumber\\
&&\IEEEeqnarraymulticol{3}{l}{\text{(\ref{eq:ps_p_constraints}) - (\ref{eq:prob_r1_constr})},} \IEEEyessubnumber
\end{IEEEeqnarray} %
and ${C}^{elec}_d$ as in (\ref{eq:c_elec}).
The peak power $P^{peak}_d$ after day $d$ is the maximum of $\bm{P}^{grid}_d$, the grid power of day $d$, and $P^{peak}_{d-1}$. 
The expected value operator in (\ref{eq:DP_obj}) can be approximated using the SAA (\ref{eq:stoch_opt}) and the scenario reduction techniques explained in section~\ref{sec:stoch_opt}.
The final value function $V_{n_d+1}$ used in the objective of day $n_d$, the last day of the billing period is:
\begin{equation}\label{eq:V_lastday}
V_{n_d+1}(P^{peak}_{n_d}) = c_{peak} P^{peak}_{n_d}.
\end{equation}
%


With the final value function $V_{n_d+1}$ defined, we can calculate $V_{d}(P^{peak}_{d-1})$ for each day $d$ by solving (\ref{eq:DP}) recursively.
However, this value function would assume the perfect hindsight solution of the second stage peak shaving problem (see section \ref{sec:pscontrol}) and not take into account the suboptimality of a practical, non-anticipative controller.
Therefore, when solving (\ref{eq:DP}), we will instead use $V^{rule}_{d+1}$, the value of the objective (\ref{eq:DP_obj}) at day $d+1$ evaluated using the rule-based peak shaving controller of Algorithm~\ref{alg:rulebased}.

All elements of the dynamic programming control scheme are combined in Figure~\ref{fig:overview}.
At the start of day $d$, the peak power $P^{peak}_{d-1}$ is known and used as an input into the stochastic optimisation (\ref{eq:DP}), which uses $\tilde{V}^{rule}_{d+1}(P^{peak}_{d})$, a convex approximation of the value function of the next day, evaluated with the rule-based controller. 
Solving (\ref{eq:DP}) gives the FCR capacity $r_d$ and recharge controller $D_d$, used in the FCR recharge controller (\ref{eq:state_feedback}), and the peak shaving boundaries $(\bm{E}_{min,d}^{ps}, \bm{E}_{max,d}^{ps}, \bm{P}_{min,d}^{ps}, \bm{P}_{max,d}^{ps})$ from (\ref{eq:FCR_ps_div}) used in the peak shaving controller of Algorithm~\ref{alg:rulebased}. 

\subsubsection{Value Function Approximation}

\begin{figure}
	\centering
	\includegraphics[width=\columnwidth]{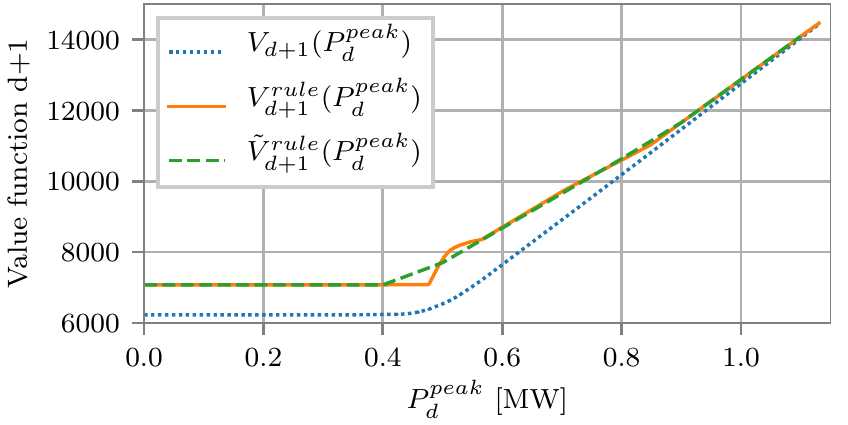}
	\caption{Example of $V_{d+1}(P^{peak}_{d})$, the value function of the optimisation (\ref{eq:DP}), $V^{rule}_{d+1}(P^{peak}_{d})$, the value function with the rule-based controller, and $\tilde{V}^{rule}_{d+1}(P^{peak}_{d})$, the convex piecewise linear approximation hereof.}
	\label{fig:valuefuncs}
\end{figure}

To solve (\ref{eq:DP}) efficiently, we need a representation of the value function $V_{d+1}^{rule}(P^{peak}_{d})$ that does not jeopardise the tractability of the optimisation problem.
As the minimisation in (\ref{eq:DP}) is convex, the value function $V_{d+1}$ is a convex function of $P^{peak}_{d}$. 
However, the value function $V_{d+1}^{rule}(P^{peak}_{d})$ is not necessarily convex, as shown in Figure \ref{fig:valuefuncs}, owing to the non-convex peak shaving controller from Algorithm \ref{alg:rulebased}. 
Therefore, we approximate $V_{d+1}^{rule}$ by a convex piecewise linear function $\tilde{V}_{d+1}^{rule}(P^{peak}_{d})$ using a least squares fit over the range $\left[0, \max(\bm{P}_{prob})\right]$.
An example of $\tilde{V}^{rule}_{d+1}$ is also shown in Figure \ref{fig:valuefuncs}. It is also interesting to note that the difference between  $\tilde{V}^{rule}_{d+1}$, $V_{d+1}^{rule}$ and $V_{d+1}$ disappears with higher $P^{peak}_{d}$, meaning that the rule based controller achieves perfect hindsight performance when $P^{peak}_{d}$ is high.

\section{Extension Towards Multiple Sites}\label{sec:mult_sites}
When multiple batteries are installed at different sites of which the shape of the consumption profiles are complementary, there can be added value in aggregating their frequency control capacity. 
For example, if one site has a high consumption peak in the morning and another site in the afternoon, the battery at the first site can do peak shaving in the morning while the battery at the second site delivers the frequency control capacity, and vice versa in the evening.

The framework for peak shaving and frequency control proposed in section \ref{sec:prob} can easily be extended to incorporate multiple sites. 
As peak tariffs are charged to each site separately, the peak shaving objective for multiple sites is simply the sum of the peak shaving objectives of the individual sites: $\min \sum_i^{n_s} c_{peak} P^{peak}_{i} + C^{elec}_i $,  with $n_s$ the number of sites. 
With regard to frequency control, the aggregated FCR capacity of all sites $r$ can be split into $n_s$ FCR capacity vectors $\bm{r}_i = (r_{1,i}, \ldots, r_{n_t,i})^T$, $i=1,\ldots, n_s$, so that the local FCR capacity can vary over time. Each site will be also have its individual recharging controller $D_i$. Finally, the individual FCR capacities have to add up to the aggregated FCR capacity $r$ at every time step $k$:
\begin{equation}\label{eq:split_r1}
r = \sum_i^{n_s} r_{k,i}, \quad \forall k = 1, \ldots, n_t.
\end{equation}
%

The dynamic programming-based control scheme of Figure~\ref{fig:overview} can  also be used for multiple site. 
Because the problem is linked by (\ref{eq:split_r1}), the value function of day $d$ is a function of $P^{peak}_{i, d-1}, i=1, \ldots, n_s$,  the peak power after day $d-1$ of every site $i$. The stochastic optimisation of (\ref{eq:DP}) becomes then:
\begin{IEEEeqnarray}{ll,rCl}
\IEEEeqnarraymulticol{5}{l}{V_d(P^{peak}_{0, d-1}, \ldots, P^{peak}_{n_s, d-1})}\nonumber\\* \;
= \> &\min& \IEEEeqnarraymulticol{3}{l}{ \E[\tilde{V}^{rule}_{d+1}(P^{peak}_{0, d}, \ldots P^{peak}_{n_s, d}) + \sum_{i=1}^{n_s} {C}^{elec}_{i,d}] -c^{FCR}_d r_d ,} \IEEEnonumber\\* \IEEEeqnarraynumspace \IEEEnonumber \label{eq:multsites}  \label{eq:multsites_obj} \\
&\mathrm{s.t.}&P^{peak}_{i,d} &\geq& {P}^{prof}_{k,i,d} + {P}^{ps}_{k,i,d} + \sum_l^{k-1} d_{kl}^{i,d} {\Delta f}_{l} + r_{k,i,d} {\Delta f}_{k}, \IEEEnonumber\\
&&&& \quad k = 1,\ldots,n_k, \quad i = 1,\ldots, n_s, \IEEEyesnumber\\
&&P^{peak}_{i,d} &\geq& P^{peak}_{i, d-1}, \quad i = 1,\ldots, n_s, \IEEEnonumber\\
&&r_d&=&\sum_{i=1}^{n_s} r_{k,i,d} , \quad k = 1, \ldots, n_k,  \IEEEnonumber\\
&&\IEEEeqnarraymulticol{3}{l}{\text{(\ref{eq:ps_p_constraints}) - (\ref{eq:prob_r1_constr})},\quad i = 1,\ldots, n_s,} \IEEEnonumber
\end{IEEEeqnarray} 

As the dimension of the state  $(P^{peak}_{0, d-1}, \ldots, P^{peak}_{n_s, d-1})$ of the dynamic program (\ref{eq:multsites}) is equal to the number of sites to be aggregated, the computational effort needed to solve the dynamic program increases with the number of sites considered~\cite{BertsekasADP}.
This can partly be mitigated by solving (\ref{eq:multsites}) for multiple states in parallel. 
More efficient sampling of the value function, using Latin hypercube sampling~\cite{cervellera2007neural} or orthogonal arrays~\cite{chen1999application} can further reduce the required computational effort when aggregating a larger number of sites, and interesting future work consists of analysing which of these methods show the best performance for the proposed problem.


\section{Simulation and Results}\label{sec:simres}
In this section we present a case study, applying the previously presented methodology to two \SI{1}{MW}, \SI{1}{MWh} batteries at two industrial sites: a pumping station (site 1) and a cold store (site 2), to perform peak shaving at the sites while delivering an aggregated FCR capacity. 
We use real consumption data from actual industrial sites and real grid frequency measurements from the CE synchronous area. 
The $5^{\text{th}}$ and $95^{\text{th}}$ percentiles of the consumption profiles are depicted by the grey shades in Figures \ref{fig:Boundaries}c and \ref{fig:Boundaries}d. The average profiles are also shown. 
The profiles are somewhat complementary: site~1 has a high peak around 7am and some lower peaks in the day 
while site~2 has the highest consumption overnight. 
The size of the batteries is chosen to be of the same order of magnitude as the consumption profiles, however, the results can be easily scaled towards sites with a larger consumption and larger batteries.

We assume the efficiencies at $\eta^c=\eta^d=\sqrt{\SI{90}{\%}}$. 
We discretise each day into time steps of 15 minutes, so $\Delta t =\SI{900}{s}$ and $n_t=96$.
In the second-order cone constraint (\ref{eq:R1_constr}), we set $\epsilon=5\cdot10^{-3}$ and calculate $\sigma_{f_k}$ and $\sigma_{b_k}$ using four year of CE frequency data. In the stochastic optimisation (\ref{eq:DP}), we draw 1500 iid scenarios which we reduce to $n_{v_r} = n_{w_r} = 50$ to obtain an SAA optimality gap $< \SI{2.5}{\%}$, following Figure \ref{fig:scen_red_SAAgap}.
We set $c^{FCR} = \SI{12}{\EUR/MW/h}, c_{peak} = \SI{13000}{\EUR/MW_{peak}/month}$ and $c_{elec}= \SI{45}{\EUR/MWh}$. 

\subsection{Combining Peak Shaving and Frequency Control}

\begin{figure}
	\centering
	\includegraphics[width=\columnwidth]{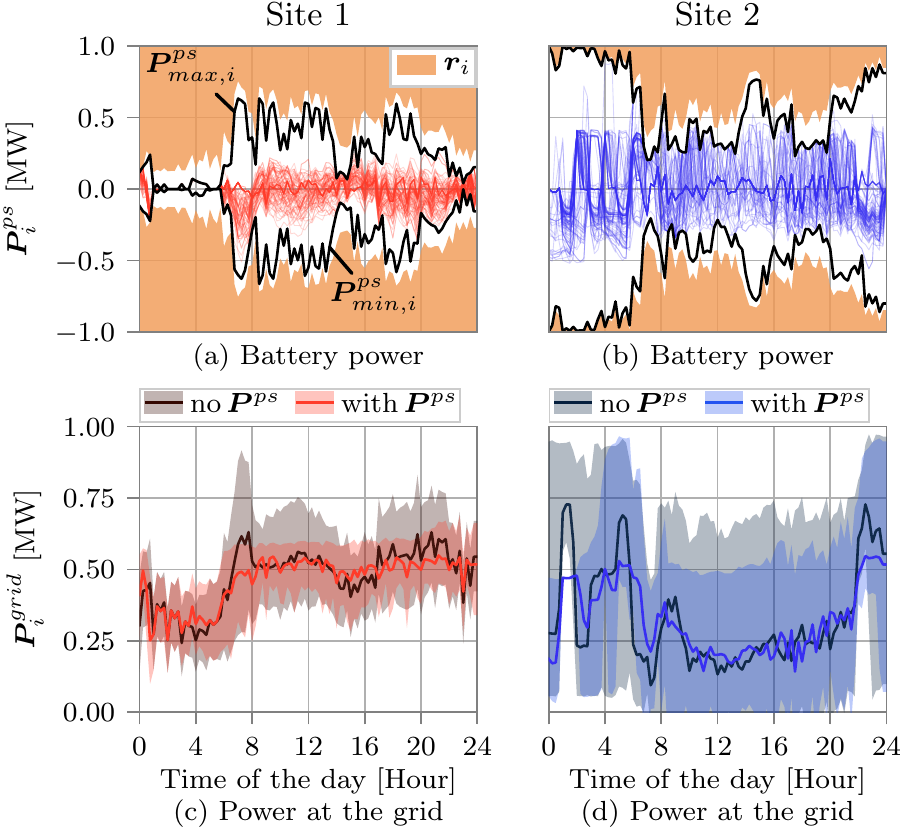}
	\caption{Aggregating two batteries at two sites performing peak shaving and frequency control. Each line in (a) \& (b) shows the battery peak shaving power for a specific scenario, while (c) \& (d) show the average (solid lines) and the $5^{\text{th}}$ - $95^{\text{th}}$ percentiles (shaded areas) of the grid power with and without the peak shaving from (a) \& (b). The coloured area in (a) \& (b) is the FCR capacity $\bm{r}_i$ of each site, which aggregated forms a constant capacity of \SI{0.88}{MW}.}
	\label{fig:Boundaries}
\end{figure}

Figure \ref{fig:Boundaries} shows how the stochastic optimisation  (\ref{eq:multsites}) succeeds in aggregating FCR capacity of the two batteries while performing peak shaving at the two sites. 
The two coloured areas in Figures \ref{fig:Boundaries}a and \ref{fig:Boundaries}b represent the FCR capacity of the sites $\bm{r}_i$, which add up to form a constant aggregated FCR capacity $r=\SI{0.88}{MW}$. 
However, at times when consumption at site 1 is expected to be high, mainly during the day, this battery delivers less FCR capacity and has more power for peak shaving available while at site~2, which has a higher consumption at night, one can see the opposite behaviour.


The coloured lines in Figures \ref{fig:Boundaries}a and \ref{fig:Boundaries}b show the actual peak shaving power scenarios $\bm{P}^{ps}_{i}$ for different daily consumption profiles $\bm{P}^{prof}_{i}$ when using the rule-based peak shaving controller of Algorithm \ref{alg:rulebased}. 
The effect of this peak shaving power on the original profiles is depicted by the coloured profiles of Figures \ref{fig:Boundaries}c and \ref{fig:Boundaries}d.
It is clear that the peak shaving power of the battery at site 1 is able to decrease the peak consumption. 
At site 2 it is more difficult to reduce the peak, as the energy content needed during to shave the peak in the first hours of the day can be more than the energy content of the battery. 
This explains the peak of the $95^{\text{th}}$ percentile around 5am-6am.
Nevertheless, the averaged profile with peak shaving 
is lower during these hours, indicating that in many scenarios the consumption power can still be reduced.



\subsection{Dynamic Programming Framework}

\begin{figure}
	\centering
	\subfloat{\centering \hspace{0.6cm}
		\includegraphics[width=0.90\columnwidth]{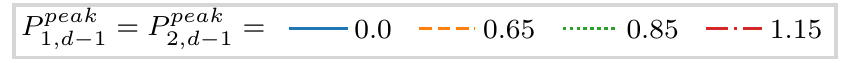}
	}\vspace*{-0.3cm}
	\setcounter{subfigure}{0}
	\subfloat[Value function $V_d$ of the optimisation (\ref{eq:multsites}) and $V_d^{rule}$ of the rule-based controller (Alg. \ref{alg:rulebased}). \label{fig:value_month}]{\centering
		\includegraphics[width=0.47\columnwidth]{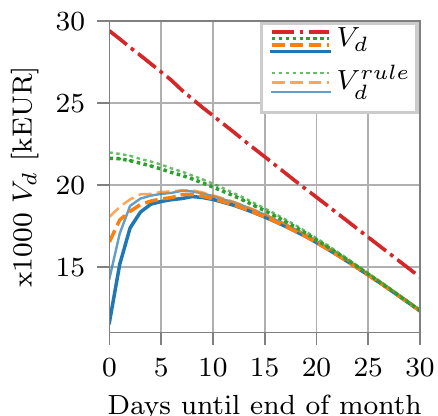}
	}
	\hfill
	\subfloat[Combined FCR capacity $r_d$ of the two batteries at the two sites. \label{fig:fcr_cap_month}]{%
		\includegraphics[width=0.47\columnwidth]{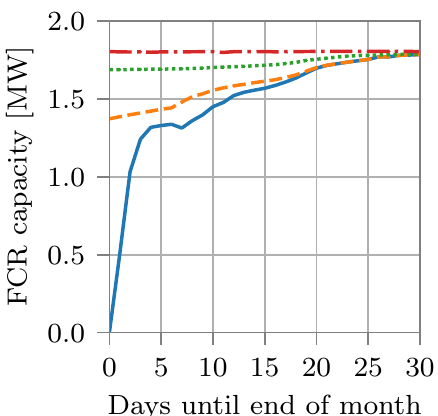}
	}
	\caption{Evolution of (a) the value function and (b) FCR capacity $r_d$ in function of the days left until the end of the month, for various peak powers $P^{peak}_{i,d-1}$ observed up to that day (here 
		equal at the two sites: $P^{peak}_{1,d-1} = P^{peak}_{2, d-1}$).}
	\label{fig:fcr_v_month}
\end{figure}

We will now look at the evolution of the decision making over time in the dynamic programming scheme.
Figure \ref{fig:value_month} shows the evolution of the value function of the dynamic program (\ref{eq:multsites}) applied to the two sites, from the last day to the first day of the month. 
The figure shows this evolution for various values of $P^{peak}_{i,d-1}$, the sum of the maximum power observed so far at the two sites. 
Figure \ref{fig:fcr_cap_month} shows the evolution of the corresponding aggregated FCR capacity $r_d$, also from the last to the first day of the month.

Analysing both figures, we can draw some insightful conclusions. 
From Figure \ref{fig:fcr_cap_month}, it turns out that a higher value of $P^{peak}_{i,d-1}$ results in a higher FCR capacity.
In case a high value of $P^{peak}_{i,d-1}$ has been observed, there is a low probability that the consumption profile will be even higher and therefore, a larger share of the battery will be allocated for FCR.
At a very high power peak $P^{peak}_{i,d-1}=1.15$, the batteries will provide their maximum FCR capacity (\SI{1.80}{MW}) over the entire month.
The linear decrease of the value function $V_d$ in Figure~\ref{fig:value_month} is thus solely due to the accumulation of FCR revenues.

\begin{figure}
	\centering
	\includegraphics[width=\columnwidth]{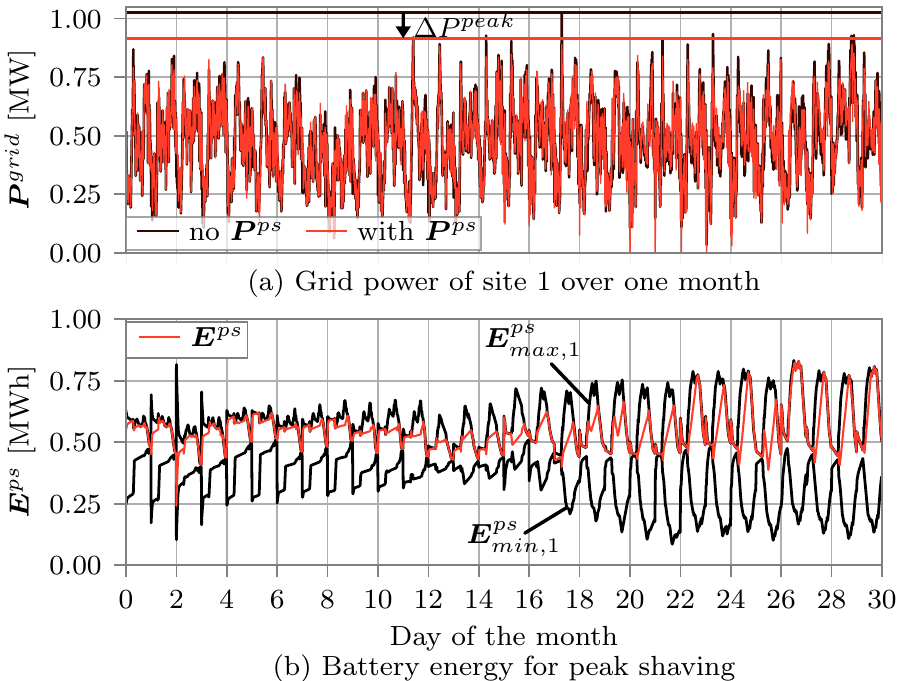}
	\caption{Example of peak shaving during one month at site 1. 
		Figure (a) shows that, with the peak shaving power of the battery, the site is able to reduce its peak power with \SI{110}{kW}.
		Figure (b) shows the corresponding actual $\bm{E}^{ps}$ and available $\bm{E}_{max}^{ps} - \bm{E}_{min}^{ps}$ energy content of the battery for peak shaving.}
	\label{fig:monthlypeak}
\end{figure}

Even in case $P^{peak}_{i,d-1}$ is low, the FCR capacity increases when more days are remaining until the end of the peak shaving period (one month). 
The longer the remaining period, the higher the probability on a high consumption peak which cannot be shaved successfully by the battery. Therefore, it is better not to lose the potential value from FCR and already use a major part of the battery for FCR. The value function of a low $P^{peak}_{i,d-1}$ will then decrease due to the FCR revenues, at almost the same rate as the value function of a high $P^{peak}_{i,d-1}$.
When the remaining period shortens and $P^{peak}_{i,d-1}$ has been rather low, there is less probability a high peak will occur in the remaining period, and the FCR capacity will be reduced as a larger share of the battery will be assigned for peak shaving trying to maintaining $P^{peak}_{i,d-1}$ low. The value function decreases as the probability of a low peak over the entire month increases.
Figure \ref{fig:monthlypeak}b, showing the peak shaving boundaries $(E_{max, ps}, E_{min, ps})$ over an entire month corresponding to the consumption profile of Figure \ref{fig:monthlypeak}a, also depicts this evolution. From the second half of the month, less capacity is used for FCR while the available energy for peak shaving becomes larger, trying to maintain the peak consumption at the level seen so far. 


Finally, we note that the value functions $V^{rule}_d$ of the rule-base peak shaving controller in Figure \ref{fig:value_month} are very close to the value functions $V^d$ from the optimisation. Except when there are few days remaining and $P^{peak}_{i,d-1}$ is low, a state which does not occur in practice, the difference becomes larger. 
Hence, we can conclude that using perfect hindsight in the second stage of the stochastic optimisation is in practice a good approximation to the actual rule-based peak shaving controller.

\subsection{Monthly Costs and Revenues}
The performance of the entire control scheme can be evaluated in Table \ref{tbl:costs}, which compares various costs components of the two sites for various cases: without batteries, with batteries performing peak shaving only, batteries performing FCR only and batteries combining FCR with peak shaving. 
The table gives the average of various consumption and frequency scenarios of one month.
The ``Peak Power" column gives the expected peak consumption power during one month of the two sites combined, while the ``Average FCR Capacity" shows the averaged FCR capacity over the entire month 
and the ``Total Net Profits" is the difference of the peak costs in the ``Without Batteries" scenario and the peak and electricity costs minus the FCR revenues of the other scenarios.

\begin{table}[ht]
	\vspace{-0.3cm}
	\centering
	\renewcommand{\arraystretch}{1.4}
	\caption{Expected Monthly Costs and Revenues of the Two Sites \newline With and Without Batteries.}
	\label{tbl:costs}
	\newcommand{\hs}{\hspace{0.0085\columnwidth}}
	\begin{tabular}{@{\hs}L{0.18\columnwidth}@{\hs}*{2}{@{\hs}C{0.09\columnwidth}@{\hs}}@{\hs}C{0.175\columnwidth}@{\hs}@{\hs}C{0.13\columnwidth}@{\hs}@{\hs}C{0.07\columnwidth}@{\hs}@{\hs}@{\hs}C{0.13\columnwidth}@{\hs}}
		\toprule
		Scenario  & Peak Power [MW] & Peak Costs [k\EUR] &  Average FCR Capacity [MW] &  FCR Revenues [k\EUR] & Elec. Costs [\EUR] & Total Net Profits [k\EUR] \\ \midrule
		Without Batteries & 1.91 & 24.9 & -- & -- & -- & -- \\ 
		Only Peak Shaving  & 1.35 & 17.5 & -- & -- & 197 & 7.2 \\ 
		Only FCR & 2.09 & 27.3 & 1.80 & 15.6 & 118 & 13.1  \\	
		\textbf{FCR \& Peak Shaving Combined} &\textbf{1.96} & \textbf{25.5} & \textbf{1.76} & \textbf{15.2} & \textbf{177} & \textbf{14.4} \\
		\bottomrule
	\end{tabular}
\end{table}


When only performing peak shaving, the batteries are able to reduce the power peak with \SI{560}{kW}, which results in a decrease of peak power costs of \SI{7200}{\EUR}. 
When only performing frequency control, the batteries together provide the maximum FCR capacity of \SI{1.80}{MW} during the entire month, which gives a revenue of \SI{15600}{\EUR}. However, this also leads to an increase in peak power to \SI{2.09}{MW}, reducing the net profits to \SI{13100}{\EUR}.
However, when combining FCR and peak shaving using the proposed methodology, the batteries are able to maintain the peak power at \SI{1.96}{MW}, while still providing \SI{1.76}{MW} of FCR capacity on average, resulting in a net profit of \SI{14400}{\EUR/month}.
In all scenarios, the additional electricity costs $C^{elec}$ of the batteries are negligible.

\section{Conclusion}\label{sec:conclusion}

In this paper, we have proposed a novel stochastic optimisation and control framework that is able to optimally combine peak shaving and frequency control objectives with a battery system installed behind-the-meter.
The framework also allows to aggregate frequency control capacity of multiple batteries at different sites, thereby leveraging potential synergies.

In a case study on two \SI{1}{MW}, \SI{1}{MWh} batteries at two industrial sites, we show that combining peak shaving with frequency control using the proposed optimisation framework leads to an expected monthly profit of \SI{14400}{\EUR}, which is two times the profit in case they would only perform peak shaving and around \SI{10}{\%} more than only performing frequency control.  





%
\bibliographystyle{IEEEtran}
\bibliography{IEEEabrv,bibl_r1_peakshaving}{}

%
%
\begin{IEEEbiography}[{\includegraphics[width=1in,height=1.25in,clip,keepaspectratio]{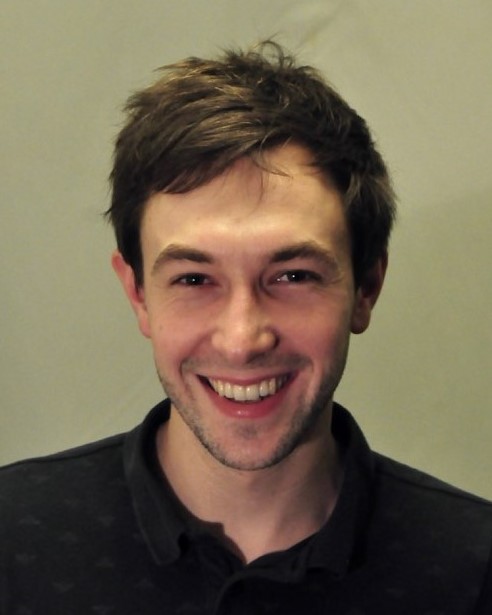}}]{Jonas Engels} (S'16), born in Lier, Belgium, received 
	his M.Sc degree in Energy Engineering from KU Leuven, Belgium in 2014.
	Afterwards, he started working at a consultancy where he worked for various major utilities.
	Currently, he works as a researcher at Centrica Business Solutions Belgium while pursuing a Ph.D. degree at the research group ELECTA of the department of Electrical Engineering (ESAT) at KU Leuven, Belgium. 
	His research is funded by the Flanders Innovation \& Entrepreneurship agency (VLAIO).
	His main research interests are in control algorithms and market optimization for smart grids, demand response and battery storage systems.
\end{IEEEbiography}

\begin{IEEEbiography}[{\includegraphics[width=1in,height=1.25in,clip,keepaspectratio]{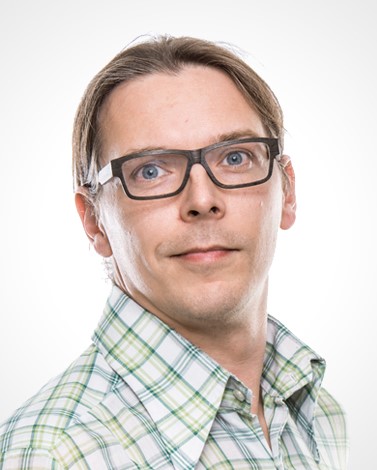}}]{Bert Claessens}
	received his M.Sc. and Ph.D. degrees in applied	physics from the Eindhoven University of Technology, The Netherlands, in 2002 and 2006, respectively.
	In 2006 he started working at ASML Veldhoven, the Netherlands, as a design engineer. Starting in 2010, he worked six years as a Researcher at the Vlaamse Instelling voor Technologisch Onderzoek	(VITO), Mol, Belgium. 
	Currently, he is head of research at Centrica Business Solutions Belgium and Full Professor at the Eindhoven University of Technology. 
	His research main interests are directed towards residential demand response and applying state of the art in artificial intelligence for energy applications.
\end{IEEEbiography}

\vfill

\newpage
%

\begin{IEEEbiography}[{\includegraphics[width=1in,height=1.25in,clip,keepaspectratio]{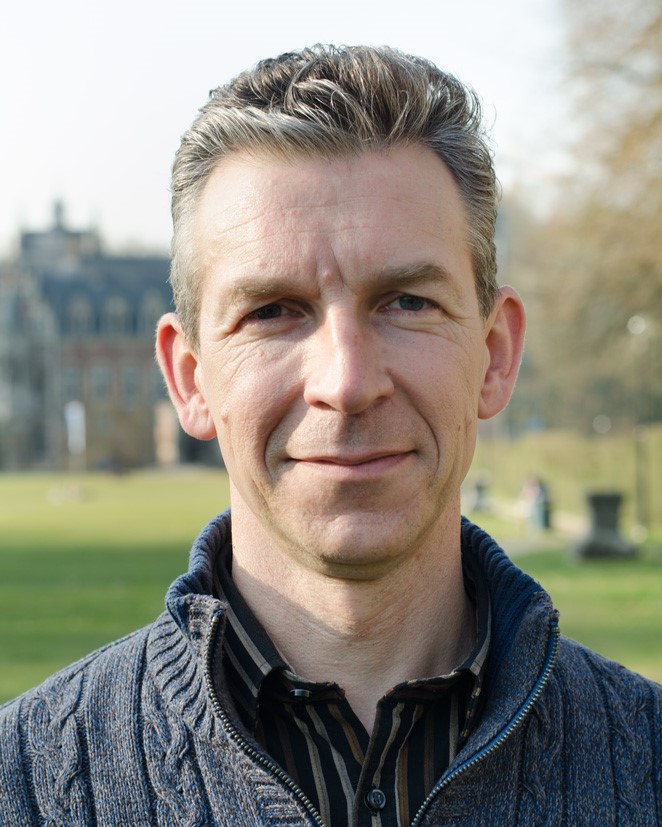}}]{Geert Deconinck} (SM'00) is full professor (gewoon hoogleraar) at KU Leuven, Belgium. 
	He received his M.Sc. degree in Electrical Engineering and his Ph.D. degree in Engineering Sciences from KU Leuven, Belgium in 1991 and 1996 respectively.
	
	He is head of the research group ELECTA on Electrical Energy at the Department of Electrical Engineering (ESAT). 
	In the research centre EnergyVille on smart energy for sustainable cities, he is the scientific leader for the ‘algorithms, modelling, optimisation’, applied to smart electrical and thermal networks. His research focuses on robust distributed coordination and control, specifically in the context of smart grids. 
	He is a fellow of the Institute of Engineering and Technology (IET).
\end{IEEEbiography}


\vfill

\end{document}